\documentclass[11pt]{article}
\usepackage{amssymb}
\usepackage{epsf}

\def\binom#1#2{{#1\choose #2}}
\newcommand{\vc}[1]{\nca{#1}}

\newcommand{\vcg}[1]{\mbox{\boldmath $#1$}}
\newcommand{\vcgs}[1]{\mbox{\scriptsize\boldmath $#1$}}

\newcommand{\nca}[1]{\mbox{\bf #1}}
\newcommand{\icas}[1]{\mbox{\scriptsize \it #1}}

\newcommand{\A}{\mathbb{A}}
\newcommand{\C}{\mathbb{C}}
\newcommand{\F}{\mathbb{F}}
\newcommand{\hache}{\mathbb{H}}

\newcommand{\K}{\mathbb{K}}
\newcommand{\N}{\mathbb{N}}
\newcommand{\og}{\mathbb{O}}

\newcommand{\R}{\mathbb{R}}
\newcommand{\Z}{\mathbb{Z}}

\title{
Basic Calculations on Clifford Algebras 
}
\author{Guillermo Morales-Luna \\
Computer Science Department \\
CINVESTAV-IPN, Mexico \\
{\tt gmorales@cs.cinvestav.mx} }

\begin{document}
\maketitle

\begin{abstract}
Clifford algebras are important structures in Geometric Algebra and Quantum Mechanics. They have allowed a formalization of the primitive operators in Quantum Theory. The algebras are built over vector spaces with dimension a power of 2 with addition and multiplication being effectively computable relative to the computability of their own spaces. Here we emphasize the algorithmic notions of the Clifford algebras. We recall the reduction of Clifford algebras into isomorphic structures also suitable for symbolic manipulation.
\end{abstract}

\section{Notation}

Let $\K=\R,\C$ be the field of real or complex numbers. $\K^n$ stands for the $n$-Cartesian power of $\K$, and $\K^{n\times n}$ for the space of square $(n\times n)$-matrices with entries in $\K$. Both $\K^n$ and $\K^{n\times n}$ are vector spaces over $\K$ with natural structures. $\K^{n\times n}$ is an algebra. With respect to matrix multiplication, the following collections of matrices are subgroups of general use:
\begin{description}
\item[General Linear Group]  $\mbox{GL}(n,\K) = \{A\in\K^{n\times n}|\ \det(A)\not=0\}$.
\item[Orthogonal Group]  $\mbox{O}(n,\K) = \{A\in\mbox{GL}(n,\K)|\ A^HA=\vc{1}\}$.
\item[Special Orthogonal Group]  $\mbox{SO}(n,\K) = \{A\in\mbox{O}(n,\K)|\ \det(A)=1\}$.
\end{description}
Some other particular groups are the following:
\begin{itemize}
\item  $\mbox{O}(n) = \mbox{O}(n,\R)$  
\item  $\mbox{SO}(n) = \mbox{SO}(n,\R)$  
\item  $\mbox{U}(n) = \mbox{O}(n,\C)$  
\item  $\mbox{SU}(n) = \mbox{SO}(n,\C)$  
\end{itemize}
For any $m,n\in\N$, by $[\![m,n]\!]$ we will denote the set of integers $\{m,m+1,\ldots,n-1,n\}$.

\section{Division Algebras}

Let $\R$ be the field of real numbers. Let $i=\sqrt{-1}$ be a square root of $-1$ and let $\C=\R[i]$ be the field of complex numbers.

\subsection{Quaternions}

Let $i,j,k$ be three symbols with relations
\begin{equation}
\begin{array}{cc}
 \multicolumn{2}{c}{i^2 = j^2 = k^2 = ijk = -1}  \\
\begin{array}{rcl}
ij &=& k \\
jk &=& i \\
ki &=& j
\end{array} & 
\begin{array}{rcl}
ji &=& - k \\
ik &=& - j \\
kj &=& - i
\end{array} 
\end{array} \label{eq.gn12}
\end{equation}
The {\em quaternion algebra} is $\hache = \R[i,j,k]$. It is a non-commutative associative division algebra, extending $\C$, and it is a 4-dimensional real vector space. 

The {\em conjugate} map is $\hache\to\hache$, $h = x_0 + x_1 i + x_2 j + x_3 k \mapsto \overline{h} = x_0 - x_1 i - x_2 j - x_3 k$ and it is congruent with respect to addition and multiplication. The {\em norm} of a quaternion $h\in\hache$ is $|h| = \sqrt{h \overline{h}}$ and, whenever it is non-zero, its multiplicative inverse is $h^{-1} = \frac{1}{|h|^2}\overline{h}$.

The {\em unit sphere} in $\hache$ is
\begin{equation}
S_{\hache} = \{h\in\hache|\ h \overline{h} = 1 \}, \label{eq.gn}
\end{equation}
and it is a subgroup under multiplication. Let $S_2 = \{x\in\R^3|\ \|x\|=1\}$ be the unit sphere in the 3-dimensional real space, and let $\phi:S_2\to S_{\hache}$, $(w_1,w_2,w_3)\mapsto w_1i+w_2j+w_3k$. Then each point $h\in S_{\hache}$ can be written in the form
\begin{equation}
h = e^{a\,\phi(w)} = \cos a + \phi(w) \sin a,\ \ \mbox{ with } -\pi < a \leq \pi\ ,\ w\in S_2. \label{eq.gn17}
\end{equation}
Let $R_3=\{x_0 + x_1 i + x_2 j + x_3 k\in\hache|\ x_0=0\}$ be the copy of $\R^3$ consisting of quaternions with zero real part. Let $A: S_{\hache}\times R_3 \to R_3$, $(h,x)\mapsto h x \overline{h}$. It is an action of the group $S_{\hache}$ over $R_3$ and for each $h = e^{a\,\phi(w)}\in S_{\hache}$, the map $x\mapsto A(h,x)$ is a counterclockwise rotation of angle $2a$ of $R_3$ along the axis $w$. In this way, it is said that $\mbox{SO}(3,\R)$ is covered twice by $S_{\hache}$.

There are several matrix representation of the quaternions. Let 
\begin{equation}
\Phi_0:\hache \to \C^{2\times 2}\ \ ,\ \ x_0 + x_1 i + x_2 j + x_3 k \mapsto \left[\begin{array}{cc}
\ x_0 + x_1 i & x_2 + x_3 i \\
 -x_2 + x_3 i & x_0 - x_1 i 
\end{array}\right] \label{eq.gn14}
\end{equation}
Then $\Phi_0$ is an embedding that preserves addition and multiplication. Similarly, let 
\begin{equation}
\Phi_1:\hache \to \R^{4\times 4}\ \ ,\ \ x_0 + x_1 i + x_2 j + x_3 k \mapsto \left[\begin{array}{rrrr}
 x_0 & -x_1 &  x_3 & -x_2 \\
 x_1 &  x_0 & -x_2 & -x_3 \\
-x_3 &  x_2 &  x_0 & -x_1 \\
 x_2 &  x_3 &  x_1 &  x_0
\end{array}\right] \label{eq.gn15}
\end{equation}
Then $\Phi_1$ is also an embedding that preserves addition and multiplication.

Hence $\hache$ is a subalgebra of both $\C^{2\times 2}$ and $\R^{4\times 4}$.

\subsection{Octonions}

The set of {\em octonions} is $\og=\hache\times\hache$ with addition defined component-wise and multiplication given by the {\em Cayley-Dickson rule}:
\begin{equation}
\left((h_0,h_1),(h_2,h_3)\right) \mapsto (h_0h_2-h_3 \overline{h_1},\overline{h_0}h_3 + h_1h_2) \label{eq.gn16}
\end{equation}
$\og$ is a non-associative division algebra, and is a real vector space of dimension 8 with basis $(1,0),(i,0),(j,0),(k,0),(0,1),(0,i),(0,j),(0,k)$. Any octonion can be written in the form $x = x_0 + \sum_{j=1}^7 x_j i_j$, with {\em conjugate} $\overline{x} = x_0 - \sum_{j=1}^7 x_j i_j$, and {\em norm} $|x|=\sqrt{x\overline{x}}$. Thus, if $x$ is non-zero, its multiplicative inverse is $x^{-1} = \frac{1}{|x|^2}\overline{x}$.

\section{Clifford Algebras}

Let $\K$ be the field $\R$ of real numbers or the field $\C$ of complex numbers. 

Let $n\in \N$, $n\geq 1$, let $\K^{2^n}$ be the $2^n$-dimensional vector space over $\K$ and let $\left(\vc{e}^{(2^n)}_j\right)_{0\leq j\leq 2^n-1}$ be its canonical basis. Each index $j\in[\![0,2^n-1]\!]$ in radix 2 is written as an $n$-length bit string: 
\begin{equation}
j = (\vcg{\varepsilon})_2 = (\varepsilon_{n-1}\cdots\varepsilon_{1}\varepsilon_{0})_2; \label{eq.gm03}
\end{equation}
and also each list corresponds to a subset $I\subset [\![0,2^n-1]\!]$: 
\begin{equation}
i\in I\ \Leftrightarrow\ \varepsilon_{i}=1. \label{eq.gm04}
\end{equation}
Thus we may write equivalently
\begin{equation}
\vc{e}^{(2^n)}_j = \vc{e}^{(2^n)}_{\vcgs{\varepsilon}} = \vc{e}^{(2^n)}_I. \label{eq.gn01}
\end{equation}
Let us denote the collection of subsets in $[\![0,n-1]\!]$ with cardinality $k$, with $0\leq k\leq n$, by 
\begin{equation}
[\![0,n-1]\!]^{(k)} =\{I\subset [\![0,n-1]\!] | \mbox{card}(I)=k\}. \label{eq.gn02}
\end{equation}
Naturally $\mbox{card}\left([\![0,n-1]\!]^{(k)}\right) = \binom{n}{k}$. 

The index set $[\![0,2^n-1]\!]$ has its canonical ordering inherited from the ordering of the integers, and this ordering corresponds to the lexicographical ordering in the set of words $\{0,1\}^n$. But on the other side, the power set ${\cal P}([\![0,n-1]\!])$ has another usual ordering $\leq_s$ defined by the following rules:
\begin{eqnarray}
I_0\in [\![0,n-1]\!]^{(k_0)}\ \&\ I_1\in [\![0,n-1]\!]^{(k_1)}\ \&\ k_0<k_1 &\Longrightarrow& I_0\leq_s  I_1 \label{eq.gm01} \\
I_0\in [\![0,n-1]\!]^{(k_0)}\ \&\ I_1\in [\![0,n-1]\!]^{(k_1)}\ \&\ k_0=k_1 &\Longrightarrow& \left[ I_0\leq_s  I_1\Leftrightarrow \vcg{\varepsilon}_0\leq_{\icas{lg}}\vcg{\varepsilon}_1\right] \label{eq.gm02}  
\end{eqnarray}
where the sets $I$ and the words $\vcg{\varepsilon}$ are related according to~(\ref{eq.gm04}) and $\leq_{\icas{lg}}$ is the lexicographical ordering. The ordering $\leq_{\icas{s}}$ can be translated into the index set $[\![0,2^n-1]\!]$ via the bijection defined by~(\ref{eq.gm03}) and~(\ref{eq.gm04}). Let $\pi$ be the permutation that sorts the canonical basis according to $\leq_{\icas{s}}$. The map $U_n:\vc{e}^{(2^n)}_j\mapsto\vc{e}^{(2^n)}_{\pi(j)}$ is unitary and for most cases its determinant is 1, i.e. it maintains the orientation of the space $\K^{2^n}$. Namely, for $n\leq 20$ the only values of $n$ for which $U_n$ changes the orientation are $3,4,5,8,9,16,17$ (integers of the form $2^k,2^k+1$ with the exception of 1 and 2).

Let us denote by
\begin{equation}
\bigwedge^k\K^n =  {\cal L}\left(\vc{e}^{(2^n)}_I\right)_{I\in[\![0,n-1]\!]^{(k)}}. \label{eq.gn03}
\end{equation}
the space spanned by the canonical vectors $\vc{e}^{(2^n)}_I$, with $I$ taken over all $k$-elements sets in $[\![0,n-1]\!]$.

$\bigwedge^1\K^n$, $\bigwedge^2\K^n$, $\bigwedge^3\K^n$ are called the spaces of {\em vectors}, {\em bivectors} and {\em trivectors} respectively. Evidently $\K^n$ is identified with $\bigwedge^1\K^n$ by the linear isomorphism such that 
\begin{equation}
\vc{e}^{(n)}_j \leftrightarrow \vc{e}^{(2^n)}_{\{j\}}. \label{eq.gn11}
\end{equation}
Let $p,q\in\N$ be such that $p+q=n$ and let $\vcg{\eta}_{pq} = \{+1\}^p \times \{-1\}^q$ be the vector $\vcg{\eta} = (\eta_0,\cdots,\eta_{n-1})$ such that 
\begin{equation}
\eta_j = \left\{\begin{array}{rl}
 1 & \mbox{ if } j<p \\
-1 & \mbox{ otherwise } 
\end{array}\right. \label{eq.gn04}
\end{equation}
$\vcg{\eta}_{pq}$ is called the {\em $(p,q)$-signature}, or {\em $(p,q)$-spacetime}, vector.

The {\em product} of the canonical vectors is defined as follows:
\begin{equation}
\begin{array}{rcl}
i<j &\Longrightarrow& \vc{e}^{(n)}_i \cdot \vc{e}^{(n)}_j =       \vc{e}^{(2^n)}_{\{i,j\}}  \\
i=j &\Longrightarrow& \vc{e}^{(n)}_i \cdot \vc{e}^{(n)}_i = \eta_i\vc{e}^{(2^n)}_{\emptyset}  \\
i>j &\Longrightarrow& \vc{e}^{(n)}_i \cdot \vc{e}^{(n)}_j = -     \vc{e}^{(n)}_j \vc{e}^{(n)}_i  
\end{array} \label{eq.gn05}
\end{equation}
and is extended to the whole basis of $\K^{2^n}$ by the following procedures
\begin{equation}
\begin{array}{rcl}
\vc{e}^{(2^n)}_{\emptyset} \cdot \vc{e}^{(n)}_j &=&  \vc{e}^{(n)}_j \\
\vc{e}^{(2^n)}_{I\cup\{i\}} \cdot \vc{e}^{(n)}_j &=&  \vc{e}^{(2^n)}_{I} \cdot \left(\vc{e}^{(n)}_i \cdot \vc{e}^{(n)}_j\right)  
\end{array} \label{eq.gn06}
\end{equation}
and
\begin{equation}
\begin{array}{rcl}
\vc{e}^{(2^n)}_{I} \cdot \vc{e}^{(2^n)}_{\emptyset} &=&  \vc{e}^{(2^n)}_{I} \\
\vc{e}^{(2^n)}_{I} \cdot \vc{e}^{(2^n)}_{J\cup\{j\}} &=& \left(\vc{e}^{(2^n)}_{I} \cdot \vc{e}^{(2^n)}_{J}\right) \cdot \vc{e}^{(n)}_j 
\end{array} \label{eq.gn07}
\end{equation}
It can be seen that the product is associative and has unit $\vc{e}^{(2^n)}_{\emptyset}=\vc{e}^{(2^n)}_{0}$.

For instance, in Table~\ref{tb.gn14} there appears the resulting product in the canonical basis of $\K^{16}$, for $n=4$ and $(p,q)=(1,3)$, and in Table~\ref{tb.gn34} there appears the resulting product, for $(p,q)=(3,1)$. In these tables at each entry $(i,j)$ a value $-k$ means $\vc{e}^{(2^n)}_i \cdot \vc{e}^{(2^n)}_j = -\vc{e}^{(2^n)}_k$ and a value $k$ means $\vc{e}^{(2^n)}_i \cdot \vc{e}^{(2^n)}_j = \vc{e}^{(2^n)}_k$. 
\begin{table}
$$\begin{array}{rrrrrrrrrrrrrrrr}
 {0} & {1} & {2} & {3} & {4} & {5} & {6} & {7} & 
 {8} & {9} & {10} & {11} & {12} & {13} & {14} & {15} \\  
 {1} & {0} & {5} & {6} & {7} & {2} & {3} & {4} & 
 {11} & {12} & {13} & {8} & {9} & {10} & {15} & {14} \\  
 {2} & - {5} & - {0} & {8} & {9} & {1} & - {11} & - {12} & 
-{3} & - {4} & {14} & {6} & {7} & - {15} & - {10} & 
 {13} \\   {3} & - {6} & - {8} & - {0} & {10} & {11} & {1} & 
-{13} & {2} & - {14} & - {4} & - {5} & {15} & {7} & {9} & 
-{12} \\   {4} & - {7} & - {9} & - {10} & - {0} & {12} & 
 {13} & {1} & {14} & {2} & {3} & - {15} & - {5} & - {6} & 
-{8} & {11} \\   {5} & - {2} & - {1} & {11} & {12} & {0} & 
-{8} & - {9} & - {6} & - {7} & {15} & {3} & {4} & - {14} & 
-{13} & {10} \\   {6} & - {3} & - {11} & - {1} & {13} & 
 {8} & {0} & - {10} & {5} & - {15} & - {7} & - {2} & {14} & 
 {4} & {12} & - {9} \\   {7} & - {4} & - {12} & - {13} & 
-{1} & {9} & {10} & {0} & {15} & {5} & {6} & - {14} & 
-{2} & - {3} & - {11} & {8} \\   {8} & {11} & {3} & - {2} & 
 {14} & {6} & - {5} & {15} & - {0} & {10} & - {9} & - {1} & 
 {13} & - {12} & - {4} & - {7} \\   {9} & {12} & {4} & 
-{14} & - {2} & {7} & - {15} & - {5} & - {10} & - {0} & 
 {8} & - {13} & - {1} & {11} & {3} & {6} \\  
 {10} & {13} & {14} & {4} & - {3} & {15} & {7} & - {6} & 
 {9} & - {8} & - {0} & {12} & - {11} & - {1} & - {2} & 
-{5} \\   {11} & {8} & {6} & - {5} & {15} & {3} & - {2} & 
 {14} & - {1} & {13} & - {12} & - {0} & {10} & - {9} & 
-{7} & - {4} \\   {12} & {9} & {7} & - {15} & - {5} & {4} & 
-{14} & - {2} & - {13} & - {1} & {11} & - {10} & - {0} & 
 {8} & {6} & {3} \\   {13} & {10} & {15} & {7} & - {6} & 
 {14} & {4} & - {3} & {12} & - {11} & - {1} & {9} & - {8} & 
-{0} & - {5} & - {2} \\   {14} & - {15} & - {10} & {9} & 
-{8} & {13} & - {12} & {11} & - {4} & {3} & - {2} & {7} & 
-{6} & {5} & {0} & - {1} \\   {15} & - {14} & - {13} & 
 {12} & - {11} & {10} & - {9} & {8} & - {7} & {6} & - {5} & 
 {4} & - {3} & {2} & {1} & - {0}
\end{array}$$
\caption{Product in the canonical basis of $\K^{16}$, for $n=4$ and $(p,q)=(1,3)$ \label{tb.gn14}}
\end{table}

\begin{table}
$$\begin{array}{rrrrrrrrrrrrrrrr}
 {0} & {1} & {2} & {3} & {4} & {5} & {6} & {7} & 
 {8} & {9} & {10} & {11} & {12} & {13} & {14} & {15} \\  
 {1} & {0} & {5} & {6} & {7} & {2} & {3} & {4} & 
 {11} & {12} & {13} & {8} & {9} & {10} & {15} & {14} \\  
 {2} & - {5} & {0} & {8} & {9} & - {1} & - {11} & - {12} & 
 {3} & {4} & {14} & - {6} & - {7} & - {15} & {10} & 
-{13} \\   {3} & - {6} & - {8} & {0} & {10} & {11} & 
-{1} & - {13} & - {2} & - {14} & {4} & {5} & {15} & - {7} & 
-{9} & {12} \\   {4} & - {7} & - {9} & - {10} & - {0} & 
 {12} & {13} & {1} & {14} & {2} & {3} & - {15} & - {5} & 
-{6} & - {8} & {11} \\   {5} & - {2} & {1} & {11} & {12} & 
-{0} & - {8} & - {9} & {6} & {7} & {15} & - {3} & - {4} & 
-{14} & {13} & - {10} \\   {6} & - {3} & - {11} & {1} & 
 {13} & {8} & - {0} & - {10} & - {5} & - {15} & {7} & {2} & 
 {14} & - {4} & - {12} & {9} \\   {7} & - {4} & - {12} & 
-{13} & - {1} & {9} & {10} & {0} & {15} & {5} & {6} & 
-{14} & - {2} & - {3} & - {11} & {8} \\  
 {8} & {11} & - {3} & {2} & {14} & - {6} & {5} & {15} & 
-{0} & - {10} & {9} & - {1} & - {13} & {12} & - {4} & 
-{7} \\   {9} & {12} & - {4} & - {14} & - {2} & - {7} & 
-{15} & - {5} & {10} & {0} & {8} & {13} & {1} & {11} & 
-{3} & - {6} \\   {10} & {13} & {14} & - {4} & - {3} & 
 {15} & - {7} & - {6} & - {9} & - {8} & {0} & - {12} & 
-{11} & {1} & {2} & {5} \\   {11} & {8} & - {6} & {5} & 
 {15} & - {3} & {2} & {14} & - {1} & - {13} & {12} & - {0} & 
-{10} & {9} & - {7} & - {4} \\   {12} & {9} & - {7} & 
-{15} & - {5} & - {4} & - {14} & - {2} & {13} & {1} & 
 {11} & {10} & {0} & {8} & - {6} & - {3} \\  
 {13} & {10} & {15} & - {7} & - {6} & {14} & - {4} & - {3} & 
-{12} & - {11} & {1} & - {9} & - {8} & {0} & {5} & {2} \\  
 {14} & - {15} & {10} & - {9} & - {8} & - {13} & {12} & 
 {11} & - {4} & - {3} & {2} & {7} & {6} & - {5} & {0} & 
-{1} \\   {15} & - {14} & {13} & - {12} & - {11} & - {10} & 
 {9} & {8} & - {7} & - {6} & {5} & {4} & {3} & - {2} & 
 {1} & - {0}
\end{array}$$
\caption{Product in the canonical basis of $\K^{16}$, for $n=4$ and $(p,q)=(3,1)$ \label{tb.gn34}}
\end{table}

Indeed, if we omit the signs in the multiplication tables over basic vectors, the resulting operation $\star : [\![0,2^n-1]\!] \times [\![0,2^n-1]\!] \to [\![0,2^n-1]\!]$ such that
\begin{equation}
i\star j = k \ \ \Longleftrightarrow \ \ \vc{e}^{(2^n)}_i \cdot \vc{e}^{(2^n)}_j = \varepsilon\, \vc{e}^{(2^n)}_k, \mbox{ for some } \varepsilon \in \{-1,1\} \label{eq.gm05}
\end{equation}
determines a structure of group over $[\![0,2^n-1]\!]$ which is isomorphic to the direct sum $\Z_2^n=\bigoplus_{i=0}^{n-1} \Z_2$ as exemplified in both tables~\ref{tb.gn14} and~\ref{tb.gn34}.

By bilinearity we extend the product over the canonical basis to a product $\K^{2^n} \times \K^{2^n} \to \K^{2^n}$:
\begin{equation}
\vc{x} = \sum_{i=0}^{2^n-1} x_i\vc{e}^{(2^n)}_i\ \&\ \vc{y} = \sum_{j=0}^{2^n-1} y_j\vc{e}^{(2^n)}_j \ \Longrightarrow\ \vc{x} \cdot \vc{y} = \sum_{i,j=0}^{2^n-1} x_iy_j\ \left(\vc{e}^{(2^n)}_i \cdot\vc{e}^{(2^n)}_j\right) = \sum_{k=0}^{2^n-1} z_k\vc{e}^{(2^n)}_k \label{eq.gn08}
\end{equation}
where the coefficients $z_k$ are determined by the coefficients $x_i$ and $y_j$. In Tables~\ref{tb.la14} and~\ref{tb.la34} we list the values of $z_k$, for $n=4$, $(p,q)=(1,3)$ and $(p,q)=(3,1)$, respectively.

\begin{table}
$$\begin{array}{c}
x_{0} y_{0} + x_{1} y_{1} - x_{2} y_{2} - x_{3} y_{3} - x_{4} y_{4} + 
x_{5} y_{5} + x_{6} y_{6} + x_{7} y_{7}  \\ 
- x_{8} y_{8} - x_{9} y_{9} - 
x_{10} y_{10} - x_{11} y_{11} - x_{12} y_{12} - x_{13} y_{13} + x_{14} y_{14} - 
x_{15} y_{15} \vspace{1ex} \\ 
x_{1} y_{0} + x_{0} y_{1} - x_{5} y_{2} - x_{6} y_{3} - 
x_{7} y_{4} + x_{2} y_{5} + x_{3} y_{6} + x_{4} y_{7}  \\ 
- x_{11} y_{8} - 
x_{12} y_{9} - x_{13} y_{10} - x_{8} y_{11} - x_{9} y_{12} - x_{10} y_{13} + 
x_{15} y_{14} - x_{14} y_{15} \vspace{1ex} \\ 
x_{2} y_{0} - x_{5} y_{1} + x_{0} y_{2} - 
x_{8} y_{3} - x_{9} y_{4} + x_{1} y_{5} - x_{11} y_{6} - x_{12} y_{7}  \\ 
+ 
x_{3} y_{8} + x_{4} y_{9} - x_{14} y_{10} - x_{6} y_{11} - x_{7} y_{12} + 
x_{15} y_{13} - x_{10} y_{14} - x_{13} y_{15} \vspace{1ex} \\ 
x_{3} y_{0} - x_{6} y_{1} + 
x_{8} y_{2} + x_{0} y_{3} - x_{10} y_{4} + x_{11} y_{5} + x_{1} y_{6} - 
x_{13} y_{7}  \\ 
- x_{2} y_{8} + x_{14} y_{9} + x_{4} y_{10} + x_{5} y_{11} - 
x_{15} y_{12} - x_{7} y_{13} + x_{9} y_{14} + x_{12} y_{15} \vspace{1ex} \\ 
x_{4} y_{0} - x_{7} y_{1} + x_{9} y_{2} + x_{10} y_{3} + x_{0} y_{4} + 
x_{12} y_{5} + x_{13} y_{6} + x_{1} y_{7}  \\ 
- x_{14} y_{8} - x_{2} y_{9} - 
x_{3} y_{10} + x_{15} y_{11} + x_{5} y_{12} + x_{6} y_{13} - x_{8} y_{14} - 
x_{11} y_{15} \vspace{1ex} \\ 
x_{5} y_{0} - x_{2} y_{1} + x_{1} y_{2} - x_{11} y_{3} - 
x_{12} y_{4} + x_{0} y_{5} - x_{8} y_{6} - x_{9} y_{7} 
 \\ 
+ x_{6} y_{8} + 
x_{7} y_{9} - x_{15} y_{10} - x_{3} y_{11} - x_{4} y_{12} + x_{14} y_{13} - 
x_{13} y_{14} - x_{10} y_{15} \vspace{1ex} \\ 
x_{6} y_{0} - x_{3} y_{1} + x_{11} y_{2} + 
x_{1} y_{3} - x_{13} y_{4} + x_{8} y_{5} + x_{0} y_{6} - x_{10} y_{7}  \\ 
- 
x_{5} y_{8} + x_{15} y_{9} + x_{7} y_{10} + x_{2} y_{11} - x_{14} y_{12} - 
x_{4} y_{13} + x_{12} y_{14} + x_{9} y_{15} \vspace{1ex} \\ 
x_{7} y_{0} - x_{4} y_{1} + 
x_{12} y_{2} + x_{13} y_{3} + x_{1} y_{4} + x_{9} y_{5} + x_{10} y_{6} + 
x_{0} y_{7}  \\ 
- x_{15} y_{8} - x_{5} y_{9} - x_{6} y_{10} + x_{14} y_{11} + 
x_{2} y_{12} + x_{3} y_{13} - x_{11} y_{14} - x_{8} y_{15} \vspace{1ex} \\ 
x_{8} y_{0} + x_{11} y_{1} - x_{3} y_{2} + x_{2} y_{3} - x_{14} y_{4} + 
x_{6} y_{5} - x_{5} y_{6} + x_{15} y_{7}  \\ 
+ x_{0} y_{8} - x_{10} y_{9} + 
x_{9} y_{10} + x_{1} y_{11} - x_{13} y_{12} + x_{12} y_{13} - x_{4} y_{14} + 
x_{7} y_{15} \vspace{1ex} \\ 
x_{9} y_{0} + x_{12} y_{1} - x_{4} y_{2} + x_{14} y_{3} + 
x_{2} y_{4} + x_{7} y_{5} - x_{15} y_{6} - x_{5} y_{7}  \\ 
+ x_{10} y_{8} + 
x_{0} y_{9} - x_{8} y_{10} + x_{13} y_{11} + x_{1} y_{12} - x_{11} y_{13} + 
x_{3} y_{14} - x_{6} y_{15} \vspace{1ex} \\ 
x_{10} y_{0} + x_{13} y_{1} - x_{14} y_{2} - 
x_{4} y_{3} + x_{3} y_{4} + x_{15} y_{5} + x_{7} y_{6} - x_{6} y_{7}  \\ 
- 
x_{9} y_{8} + x_{8} y_{9} + x_{0} y_{10} - x_{12} y_{11} + x_{11} y_{12} + 
x_{1} y_{13} - x_{2} y_{14} + x_{5} y_{15} \vspace{1ex} \\ 
x_{11} y_{0} + x_{8} y_{1} - 
x_{6} y_{2} + x_{5} y_{3} - x_{15} y_{4} + x_{3} y_{5} - x_{2} y_{6} + 
x_{14} y_{7}  \\ 
+ x_{1} y_{8} - x_{13} y_{9} + x_{12} y_{10} + x_{0} y_{11} - 
x_{10} y_{12} + x_{9} y_{13} - x_{7} y_{14} + x_{4} y_{15} \vspace{1ex} \\ 
x_{12} y_{0} + x_{9} y_{1} - x_{7} y_{2} + x_{15} y_{3} + x_{5} y_{4} + 
x_{4} y_{5} - x_{14} y_{6} - x_{2} y_{7}  \\ 
+ x_{13} y_{8} + x_{1} y_{9} - 
x_{11} y_{10} + x_{10} y_{11} + x_{0} y_{12} - x_{8} y_{13} + x_{6} y_{14} - 
x_{3} y_{15} \vspace{1ex} \\ 
x_{13} y_{0} + x_{10} y_{1} - x_{15} y_{2} - x_{7} y_{3} + 
x_{6} y_{4} + x_{14} y_{5} + x_{4} y_{6} - x_{3} y_{7}  \\ 
- x_{12} y_{8} + 
x_{11} y_{9} + x_{1} y_{10} - x_{9} y_{11} + x_{8} y_{12} + x_{0} y_{13} - 
x_{5} y_{14} + x_{2} y_{15} \vspace{1ex} \\ 
x_{14} y_{0} - x_{15} y_{1} + x_{10} y_{2} - 
x_{9} y_{3} + x_{8} y_{4} + x_{13} y_{5} - x_{12} y_{6} + x_{11} y_{7}  \\ 
+ 
x_{4} y_{8} - x_{3} y_{9} + x_{2} y_{10} - x_{7} y_{11} + x_{6} y_{12} - 
x_{5} y_{13} + x_{0} y_{14} + x_{1} y_{15} \vspace{1ex} \\ 
x_{15} y_{0} - x_{14} y_{1} + 
x_{13} y_{2} - x_{12} y_{3} + x_{11} y_{4} + x_{10} y_{5} - x_{9} y_{6} + 
x_{8} y_{7}  \\ 
+ x_{7} y_{8} - x_{6} y_{9} + x_{5} y_{10} - x_{4} y_{11} + 
x_{3} y_{12} - x_{2} y_{13} + x_{1} y_{14} + x_{0} y_{15}
\end{array}$$
\caption{Coefficients of the product in terms of the coefficients of the factors, for $n=4$ and $(p,q)=(1,3)$ \label{tb.la14}}
\end{table}

\begin{table}
$$\begin{array}{c}
x_{0} y_{0} + x_{1} y_{1} + x_{2} y_{2} + x_{3} y_{3} - x_{4} y_{4} - 
x_{5} y_{5} - x_{6} y_{6} + x_{7} y_{7}  \\ 
- x_{8} y_{8} + x_{9} y_{9} + 
x_{10} y_{10} - x_{11} y_{11} + x_{12} y_{12} + x_{13} y_{13} + x_{14} y_{14} - 
x_{15} y_{15} \vspace{1ex} \\ 
x_{1} y_{0} + x_{0} y_{1} + x_{5} y_{2} + x_{6} y_{3} - 
x_{7} y_{4} - x_{2} y_{5} - x_{3} y_{6} + x_{4} y_{7}  \\ 
- x_{11} y_{8} + 
x_{12} y_{9} + x_{13} y_{10} - x_{8} y_{11} + x_{9} y_{12} + x_{10} y_{13} + 
x_{15} y_{14} - x_{14} y_{15} \vspace{1ex} \\ 
x_{2} y_{0} - x_{5} y_{1} + x_{0} y_{2} + 
x_{8} y_{3} - x_{9} y_{4} + x_{1} y_{5} + x_{11} y_{6} - x_{12} y_{7}  \\ 
- 
x_{3} y_{8} + x_{4} y_{9} + x_{14} y_{10} + x_{6} y_{11} - x_{7} y_{12} - 
x_{15} y_{13} + x_{10} y_{14} + x_{13} y_{15} \vspace{1ex} \\ 
x_{3} y_{0} - x_{6} y_{1} - 
x_{8} y_{2} + x_{0} y_{3} - x_{10} y_{4} - x_{11} y_{5} + x_{1} y_{6} - 
x_{13} y_{7}  \\ 
+ x_{2} y_{8} - x_{14} y_{9} + x_{4} y_{10} - x_{5} y_{11} + 
x_{15} y_{12} - x_{7} y_{13} - x_{9} y_{14} - x_{12} y_{15} \vspace{1ex} \\ 
x_{4} y_{0} - x_{7} y_{1} - x_{9} y_{2} - x_{10} y_{3} + x_{0} y_{4} - 
x_{12} y_{5} - x_{13} y_{6} + x_{1} y_{7}  \\ 
- x_{14} y_{8} + x_{2} y_{9} + 
x_{3} y_{10} + x_{15} y_{11} - x_{5} y_{12} - x_{6} y_{13} - x_{8} y_{14} - 
x_{11} y_{15} \vspace{1ex} \\ 
x_{5} y_{0} - x_{2} y_{1} + x_{1} y_{2} + x_{11} y_{3} - 
x_{12} y_{4} + x_{0} y_{5} + x_{8} y_{6} - x_{9} y_{7}  \\ 
- x_{6} y_{8} + 
x_{7} y_{9} + x_{15} y_{10} + x_{3} y_{11} - x_{4} y_{12} - x_{14} y_{13} + 
x_{13} y_{14} + x_{10} y_{15} \vspace{1ex} \\ 
x_{6} y_{0} - x_{3} y_{1} - x_{11} y_{2} + 
x_{1} y_{3} - x_{13} y_{4} - x_{8} y_{5} + x_{0} y_{6} - x_{10} y_{7}  \\ 
+ 
x_{5} y_{8} - x_{15} y_{9} + x_{7} y_{10} - x_{2} y_{11} + x_{14} y_{12} - 
x_{4} y_{13} - x_{12} y_{14} - x_{9} y_{15} \vspace{1ex} \\ 
x_{7} y_{0} - x_{4} y_{1} - 
x_{12} y_{2} - x_{13} y_{3} + x_{1} y_{4} - x_{9} y_{5} - x_{10} y_{6} + 
x_{0} y_{7}  \\ 
- x_{15} y_{8} + x_{5} y_{9} + x_{6} y_{10} + x_{14} y_{11} - 
x_{2} y_{12} - x_{3} y_{13} - x_{11} y_{14} - x_{8} y_{15} \vspace{1ex} \\ 
x_{8} y_{0} + x_{11} y_{1} - x_{3} y_{2} + x_{2} y_{3} - x_{14} y_{4} + 
x_{6} y_{5} - x_{5} y_{6} + x_{15} y_{7}  \\ 
+ x_{0} y_{8} - x_{10} y_{9} + 
x_{9} y_{10} + x_{1} y_{11} - x_{13} y_{12} + x_{12} y_{13} - x_{4} y_{14} + 
x_{7} y_{15} \vspace{1ex} \\ 
x_{9} y_{0} + x_{12} y_{1} - x_{4} y_{2} - x_{14} y_{3} + 
x_{2} y_{4} + x_{7} y_{5} + x_{15} y_{6} - x_{5} y_{7}  \\ 
- x_{10} y_{8} + 
x_{0} y_{9} + x_{8} y_{10} - x_{13} y_{11} + x_{1} y_{12} + x_{11} y_{13} - 
x_{3} y_{14} + x_{6} y_{15} \vspace{1ex} \\ 
x_{10} y_{0} + x_{13} y_{1} + x_{14} y_{2} - 
x_{4} y_{3} + x_{3} y_{4} - x_{15} y_{5} + x_{7} y_{6} - x_{6} y_{7}  \\ 
+ 
x_{9} y_{8} - x_{8} y_{9} + x_{0} y_{10} + x_{12} y_{11} - x_{11} y_{12} + 
x_{1} y_{13} + x_{2} y_{14} - x_{5} y_{15} \vspace{1ex} \\ 
x_{11} y_{0} + x_{8} y_{1} - 
x_{6} y_{2} + x_{5} y_{3} - x_{15} y_{4} + x_{3} y_{5} - x_{2} y_{6} + 
x_{14} y_{7}  \\ 
+ x_{1} y_{8} - x_{13} y_{9} + x_{12} y_{10} + x_{0} y_{11} - 
x_{10} y_{12} + x_{9} y_{13} - x_{7} y_{14} + x_{4} y_{15} \vspace{1ex} \\ 
x_{12} y_{0} + x_{9} y_{1} - x_{7} y_{2} - x_{15} y_{3} + x_{5} y_{4} + 
x_{4} y_{5} + x_{14} y_{6} - x_{2} y_{7}  \\ 
- x_{13} y_{8} + x_{1} y_{9} + 
x_{11} y_{10} - x_{10} y_{11} + x_{0} y_{12} + x_{8} y_{13} - x_{6} y_{14} + 
x_{3} y_{15} \vspace{1ex} \\ 
x_{13} y_{0} + x_{10} y_{1} + x_{15} y_{2} - x_{7} y_{3} + 
x_{6} y_{4} - x_{14} y_{5} + x_{4} y_{6} - x_{3} y_{7}  \\ 
+ x_{12} y_{8} - 
x_{11} y_{9} + x_{1} y_{10} + x_{9} y_{11} - x_{8} y_{12} + x_{0} y_{13} + 
x_{5} y_{14} - x_{2} y_{15} \vspace{1ex} \\ 
x_{14} y_{0} - x_{15} y_{1} + x_{10} y_{2} - 
x_{9} y_{3} + x_{8} y_{4} + x_{13} y_{5} - x_{12} y_{6} + x_{11} y_{7}  \\ 
+ 
x_{4} y_{8} - x_{3} y_{9} + x_{2} y_{10} - x_{7} y_{11} + x_{6} y_{12} - 
x_{5} y_{13} + x_{0} y_{14} + x_{1} y_{15} \vspace{1ex} \\ 
x_{15} y_{0} - x_{14} y_{1} + 
x_{13} y_{2} - x_{12} y_{3} + x_{11} y_{4} + x_{10} y_{5} - x_{9} y_{6} + 
x_{8} y_{7}  \\ 
+ x_{7} y_{8} - x_{6} y_{9} + x_{5} y_{10} - x_{4} y_{11} + 
x_{3} y_{12} - x_{2} y_{13} + x_{1} y_{14} + x_{0} y_{15}
\end{array}$$
\caption{Coefficients of the product in terms of the coefficients of the factors, for $n=4$ and $(p,q)=(3,1)$ \label{tb.la34}}
\end{table}

For the special case of vectors $\vc{x} = \sum_{i=1}^{4} x_i\vc{e}^{(16)}_i$, $\vc{y} = \sum_{j=1}^{4} y_j\vc{e}^{(16)}_j$ $\in \bigwedge^1\K^4$ the product takes the forms listed in Table~\ref{tb.lv14}.
\begin{table}
$$\begin{array}{cc}
\begin{array}{c}
x_{1} y_{1} - x_{2} y_{2} - x_{3} y_{3} - x_{4} y_{4} \\ 0 \\ 0 \\ 0 \\ 0 \\ 
-x_{2} y_{1} + x_{1} y_{2} \\ -x_{3} y_{1} + x_{1} y_{3} \\ 
-x_{4} y_{1} + x_{1} y_{4} \\ -x_{3} y_{2} + x_{2} y_{3} \\ 
-x_{4} y_{2} + x_{2} y_{4} \\ -x_{4} y_{3} + x_{3} y_{4} \\ 0 \\ 0 \\ 0 \\ 0 \\ 0
\end{array} &
\begin{array}{c}
x_{1} y_{1} + x_{2} y_{2} + x_{3} y_{3} - x_{4} y_{4} \\ 0 \\ 0 \\ 0 \\ 0 \\ 
-x_{2} y_{1} + x_{1} y_{2} \\ -x_{3} y_{1} + x_{1} y_{3} \\ 
-x_{4} y_{1} + x_{1} y_{4} \\ -x_{3} y_{2} + x_{2} y_{3} \\ 
-x_{4} y_{2} + x_{2} y_{4} \\ -x_{4} y_{3} + x_{3} y_{4} \\ 0 \\ 0 \\ 0 \\ 0 \\ 0
\end{array} \\
\mbox{ Case: }(p,q)=(1,3) & \mbox{ Case: }(p,q)=(3,1) 
\end{array}$$
\caption{Coefficients of the product for the case of vectors, for $n=4$. \label{tb.lv14}}
\end{table}

It can be seen in general~\cite{L97} that for vectors $\vc{x} = \sum_{i=1}^{n} x_i\vc{e}^{(2^n)}_i$, $\vc{y} = \sum_{j=1}^{n} y_j\vc{e}^{(2^n)}_j$ $\in \bigwedge^1\K^n$
\begin{equation}
\frac{1}{2}\left(\vc{x}\cdot\vc{y}+\vc{y}\cdot\vc{x}\right) =  \left(\sum_{i,j=1}^{n} x_i\eta_{ij}y_j\right) \vc{e}^{(2^n)}_{\emptyset} \label{eq.gn09}
\end{equation}
where $\eta_{ij} = 0$ if $i\not=j$ and $\eta_{ii} = \eta_{i-1}$ (the sign of the $(i-1)$-entry of the $(p,q)$-spacetime vector). Besides,
\begin{equation}
\frac{1}{2}\left(\vc{x}\cdot\vc{y}-\vc{y}\cdot\vc{x}\right) =  \sum_{1\leq i<j\leq n} \frac{1}{2}\left(x_iy_j - y_ix_j\right) \vc{e}^{(2^n)}_{\{i-1,j-1\}} \label{eq.gn10}
\end{equation}

With this product, $\K^{2^n}$ is an algebra, it is called the {\em Clifford algebra} corresponding to the $(p,q)$-signature and it is denoted by $\K^{p,q}$. 

In what follows, let us recall succintly classical reduction results. The details can be found in Lounesto's book~\cite{L97}.

\section{Representation of Clifford Algebras}

\subsection{Pauli Algebra}

The Clifford algebra $\R^{3,0}$ is built over the space $\R^{2^3}$. Let $\left(\vc{e}^{(8)}_j\right)_{0\leq j\leq 7}$ be its canonical basis. The elements $\vc{e}^{(8)}_1,\vc{e}^{(8)}_2,\vc{e}^{(8)}_3$ generate the subspace of 1-vectors. The {\em Pauli matrices}
are the following on the ring $\C^{2\times 2}$:
\begin{equation}
\begin{array}{lll}
\sigma_1 = \left[\begin{array}{rr}
0 & 1 \\
1 & 0 
\end{array}\right] & 
\sigma_2 = \left[\begin{array}{rr}
0 &-i \\
i & 0 
\end{array}\right] & 
\sigma_3 = \left[\begin{array}{rr}
1 & 0 \\
0 &-1 
\end{array}\right] 
\end{array} \label{eq.gn23}
\end{equation}
and let $\Phi: \vc{e}^{(8)}_j \mapsto \sigma_j$, $j=1,2,3$. The map extends to a monomorphism $\Phi: \R^{3,0} \to \C^{2\times 2}$. Let $\sigma_0$ be the $(2\times 2)$-identity matrix.

The 4-dimensional linear subspace of $\Phi(\R^{3,0})$ with basis $\{\sigma_0 , i\sigma_1 , i\sigma_2 , i\sigma_4\}$ is a subalgebra of $\Phi(\R^{3,0})$ which is isomorphic to the quaternions by the isomorphism such that
$$\sigma_0\mapsto 1\ \ ,\ \ i\sigma_1\mapsto \hat{i}\ \ ,\ \ i\sigma_2\mapsto \hat{j}\ \ ,\ \ i\sigma_3\mapsto \hat{k},$$
where $\hat{i}$, $\hat{j}$ and $\hat{k}$ are the generator symbols of $\hache$.
Then, necessarily, 
\begin{equation}
\R^{3,0} \equiv \C\otimes\hache. \label{eq.gn24}
\end{equation}

\subsection{Minkowski's Spacetime Real Clifford Algebra}

The Minkowski's spacetime Clifford algebra $\R^{3,1}$ is built over the space $\R^{2^4}$. Let $\left(\vc{e}^{(16)}_j\right)_{0\leq j\leq 15}$ be its canonical basis. The elements $\vc{e}^{(16)}_1,\vc{e}^{(16)}_2,\vc{e}^{(16)}_3,\vc{e}^{(16)}_4$ generate the subspace of 1-vectors. Let us consider the following matrices on the ring $\R^{4\times 4}$:
\begin{equation}
\begin{array}{ll}
A_1 = \left[\begin{array}{r@{\hspace{1em}\ \ }rrr}
1 & 0 & 0 & 0 \\
0 & 1 & 0 & 0 \\
0 & 0 &-1 & 0 \\
0 & 0 & 0 &-1 
\end{array}\right] & 
A_2 = \left[\begin{array}{@{\hspace{1em}\ }r@{\hspace{1em}\ \ }rrr}
0 & 0 & 1 & 0 \\
0 & 0 & 0 & 1 \\
1 & 0 & 0 & 0 \\
0 & 1 & 0 & 0 
\end{array}\right] \vspace{2ex} \\
A_3 = \left[\begin{array}{rrr@{\hspace{1em}\ \ }r}
0 & 0 & 0 & 1 \\
0 & 0 &-1 & 0 \\
0 &-1 & 0 & 0 \\
1 & 0 & 0 & 0 
\end{array}\right] & 
A_4 = \left[\begin{array}{rrrr}
0 & 0 & 0 & 1 \\
0 & 0 & 1 & 0 \\
0 &-1 & 0 & 0 \\
-1& 0 & 0 & 0 
\end{array}\right] 
\end{array} \label{eq.gn20}
\end{equation}
and let $\Phi_R: \vc{e}^{(16)}_j \mapsto A_j$, $j=1,2,3,4$. The map extends to a monomorphism $\Phi_R: \R^{3,1} \to \R^{4\times 4}$. Thus, via $\Phi_R$, $\R^{3,1}$ is identified with a subalgebra of linear endomorphisms $\R^4\to\R^4$, and the space $\R^4$ is said to be the {\em space of Majorana spinors} of $\R^{3,1}$. 

Alternatively, let us consider the following matrices on the ring $\C^{4\times 4}$:
\begin{equation}
\begin{array}{ll}
B_1 = \left[\begin{array}{r@{\hspace{1em}\ \ }rr@{\hspace{1em}\ \ }r}
0 & 0 & 0 & 1 \\
0 & 0 & 1 & 0 \\
0 & 1 & 0 & 0 \\
1 & 0 & 0 & 0 
\end{array}\right] & 
B_2 = \left[\begin{array}{@{\hspace{1em}\ }rrrr}
0 & 0 & 0 &-i \\
0 & 0 & i & 0 \\
0 &-i & 0 & 0 \\
i & 0 & 0 & 0 
\end{array}\right] \vspace{2ex} \\
B_3 = \left[\begin{array}{rrrr}
0 & 0 & 1 & 0 \\
0 & 0 & 0 &-1 \\
1 & 0 & 0 & 0 \\
0 &-1 & 0 & 0 
\end{array}\right] & 
B_4 = \left[\begin{array}{rrr@{\hspace{1em}\ \ }r}
0 & 0 & 1 & 0 \\
0 & 0 & 0 & 1 \\
-1& 0 & 0 & 0 \\
0 &-1 & 0 & 0 
\end{array}\right] 
\end{array} \label{eq.gn21}
\end{equation}
and let $\Phi_C: \vc{e}^{(16)}_j \mapsto A_j$, $j=1,2,3,4$. The map extends to a monomorphism $\Phi_C: \R^{3,1} \to \C^{4\times 4}$. Thus, via $\Phi_C$, $\R^{3,1}$ is identified with a subalgebra of linear endomorphisms $\C^4\to\C^4$, and the space $\C^4$ is said to be the {\em space of Dirac spinors} of $\R^{3,1}$.

We may see that $A_1 = i B_1 B_2 B_3 B_4$, but, since $i$ cannot be realized as an element of $\Phi_C(\R^{3,1})$, $A_1\not\in\Phi_C(\R^{3,1})$, and the map
\begin{equation}
\Psi: \vc{e}^{(32)}_1 \mapsto B_1\ \ ,\ \ \vc{e}^{(32)}_2 \mapsto B_2\ \ ,\ \ \vc{e}^{(32)}_3 \mapsto B_3\ \ ,\ \ \vc{e}^{(32)}_4 \mapsto A_1\ \ ,\ \ \vc{e}^{(32)}_5 \mapsto B_4 \label{eq.gn22}
\end{equation}
 extends to an isomorphism $\Psi: \R^{4,1} \to \C^{4\times 4}$. 

\subsection{Opposite Metric to Minkowski's Spacetime}

The Clifford algebra $\R^{1,3}$ is built also over the space $\R^{2^4}$. Let us consider the following matrices on the ring $\hache^{2\times 2}$:
\begin{equation}
\begin{array}{llll}
\hat{I} = \left[\begin{array}{rr}
0 & \hat{i} \\
\hat{i} & 0 
\end{array}\right] & 
\hat{J} = \left[\begin{array}{rr}
0 & \hat{j} \\
\hat{j} & 0 
\end{array}\right] & 
\hat{K} = \left[\begin{array}{rr}
0 & \hat{k} \\
\hat{k} & 0 
\end{array}\right]  & 
L = \left[\begin{array}{rr}
1 & 0 \\
0 &-1 
\end{array}\right]
\end{array} \label{eq.gn25}
\end{equation}
The map $\Phi: \vc{e}^{(16)}_1 \mapsto \hat{I}$, $\vc{e}^{(16)}_2 \mapsto \hat{J}$, $\vc{e}^{(16)}_3 \mapsto \hat{K}$, $\vc{e}^{(16)}_4 \mapsto L$ extends to a monomorphism $\Phi:\R^{1,3} \to \hache^{2\times 2}$. Thus $\Phi\left(\R^{1,3}\right)$ is an algebra of linear endomorphisms $\hache^2 \to \hache^2$. The spinors here are pairs of quaternions.

We can also consider the matrices
\begin{equation}
\begin{array}{llll}
C_0 = \left[\begin{array}{rr}
0 & 1 \\
1 & 0 
\end{array}\right] & 
C_1 = \left[\begin{array}{rr}
0 &-i\hat{i} \\
i\hat{i} & 0 
\end{array}\right] & 
C_2 = \left[\begin{array}{rr}
0 &-i\hat{j} \\
i\hat{j} & 0 
\end{array}\right] & 
C_3 = \left[\begin{array}{rr}
0 &-i\hat{k} \\
i\hat{k} & 0 
\end{array}\right] 
\end{array} \label{eq.gn26}
\end{equation}
Then $L = i C_0 C_1 C_2 C_3$ and the map $\Psi: \vc{e}^{(32)}_j \mapsto C_{j-1}$, $j\in [\![1,4]\!]$, $\vc{e}^{(32)}_5 \mapsto L$, extends to a monomorphism of $\R^{2,3}$ into the complexification of $\hache^{2\times 2}$. 

As a third representation, let us consider the following matrices on $\C^{4\times 4}$:
\begin{equation}
\begin{array}{ll}
D_0 = \left[\begin{array}{rr}
0 & \sigma_0 \\
\sigma_0 & 0 
\end{array}\right] & 
D_j = \left[\begin{array}{rr}
0 &-\sigma_j \\
\sigma_j & 0 
\end{array}\right] \ \ \ j \in [\![1,3]\!]
\end{array} \label{eq.gn27}
\end{equation}
where $\sigma_j$ are the Pauli matrices.
The map $\Omega: \vc{e}^{(16)}_j \mapsto D_{j-1}$, $j \in [\![1,4]\!]$, extends to a monomorphism $\Omega:\R^{1,3} \to \C^{4\times 4}$. Thus $\Omega\left(\R^{1,3}\right)$ is an algebra of linear endomorphisms $\C^4 \to \C^4$. The spinors here are 4-tuples of complex numbers.

\subsection{Real Clifford Algebras}

Let us consider $\K=\R$. It can be seen that the following relations hold:
\begin{eqnarray}
n \equiv 0\,\mbox{\rm mod}\,2 &\Longrightarrow& \left\{\begin{array}{lcl}
(p-q) \equiv 0,2\,\mbox{\rm mod}\,8 &\Longrightarrow& \R^{p,q} \equiv \R^{2^{\frac{n}{2}}\times 2^{\frac{n}{2}}} \vspace{1ex} \\
(p-q) \equiv 4,6\,\mbox{\rm mod}\,8 &\Longrightarrow& \R^{p,q} \equiv \hache^{2^{\frac{n-2}{2}}\times 2^{\frac{n-2}{2}}} 
\end{array}\right.  \label{eq.gn18} \\
n \equiv 1\,\mbox{\rm mod}\,2 &\Longrightarrow& \left\{\begin{array}{lcl}
(p-q) \equiv 1\,\mbox{\rm mod}\,8 &\Longrightarrow& \R^{p,q} \equiv \mbox{diag}_2\left(\R^{2^{\frac{n-1}{2}}\times 2^{\frac{n-1}{2}}}\right) \vspace{1ex} \\
(p-q) \equiv 3,7\,\mbox{\rm mod}\,8 &\Longrightarrow& \R^{p,q} \equiv \C^{2^{\frac{n-1}{2}}\times 2^{\frac{n-1}{2}}} \vspace{1ex} \\
(p-q) \equiv 5\,\mbox{\rm mod}\,8 &\Longrightarrow& \R^{p,q} \equiv \mbox{diag}_2\left(\hache^{2^{\frac{n-3}{2}}\times 2^{\frac{n-3}{2}}}\right) 
\end{array}\right.   \label{eq.gn19} 
\end{eqnarray}
where, for any ring of square matrices $\A$, $\mbox{diag}_2\left(\A\right)$ consists of the matrices of the form $\left(\begin{array}{cc}
A & 0 \\
0 & B
\end{array}\right)$, with $A,B\in\A$.

\bibliographystyle{plain}
\bibliography{bascal}

\end{document}